\documentclass[a4paper,12pt]{article}

\usepackage{wrapfig}

\renewenvironment{thebibliography}[1]{
  \begin{oldthebibliography}{#1}
  	 \footnotesize
    \setlength{\itemsep}{0em}
    \setlength{\parskip}{0em}
}
{
  \end{oldthebibliography}
}

\usepackage{amsmath,amssymb,amsthm}
\numberwithin{equation}{section}

\usepackage{verbatim}
\usepackage[toc,page]{appendix}
\usepackage{graphicx}
\usepackage{fourier}
\theoremstyle{plain} 
\newtheorem{thm}{Theorem}[section] 
\newtheorem*{theorem*}{Theorem}

\theoremstyle{definition}

\newtheorem*{defn*}{Definition}

\begin{document}

\renewcommand{\title}[1]{\vspace{10mm}\noindent{\Large{\bf#1}}\vspace{8mm}} 
\newcommand{\authors}[1]{\noindent{\large #1}\vspace{5mm}} 
\newcommand{\address}[1]{{\itshape #1}}

\begin{center}
\vskip -9mm
\title{{\Large Derivation of the HOMFLYPT knot polynomial via helicity and geometric quantization}}
\vskip -3mm
\authors{Antonio Michele MITI${}^{\dagger, \ddagger}$ and Mauro SPERA${}^{\dagger}$}
\vskip -1mm
\address{
${}^{\dagger}$ Dipartimento di Matematica e Fisica  ``Niccol\`o Tartaglia", 
Universit\`a Cattolica del Sacro Cuore\\
Via dei Musei 41, 25121 Brescia, Italia\\
\vspace{2mm}
${}^{\ddagger}$ Departement Wiskunde, KU-Leuven\\ 
Celestijnenlaan 200B, B-3001 Leuven (Heverlee), Belgi\"e
}
\vskip 1mm
\textbf{Abstract}
\vskip 2mm
\begin{minipage}{15cm}%
In this Letter  we propose a semiclassical interpretation of the HOMFLYPT polynomial building on the Liu-Ricca hydrodynamical approach to the latter and on the Besana-S. symplectic approach to framing via Brylinski's manifold of mildly singular links.
\end{minipage}
\end{center}
\smallskip

\par
geometric quantization, Maslov index.

\section{Introduction}\par
In this article, building on the Maslov-type methods developed in \cite{BeSpe06}, we present a novel interpretation of the HOMFLYPT (and hence of the Jones) polynomial (\cite{Freyd-etal,PT}) as a WKB-wave function via geometric quantization of the so-called Brylinski manifold of singular knots (and links), taking inspiration from the  {\it ad hoc} helicity-based hydrodynamical procedures devised in \cite{Liu-Ricca12,Liu-Ricca15}.
Our approach can be compared with the Jeffrey-Weitsman one (\cite{Jef-Wei92,Jef-Wei93}), providing a rigorous framework for the Jones-Witten theory (\cite{Witten89,Kohno02}). The latter, though again based on geometric quantization, is much more sophisticated. In our setting, no reference to Lie groups (other than $U(1)$) is made and, as in Liu-Ricca, everything is based on helicity only, at the cost of relying on the Maslov-H\"ormander techniques of \cite{BeSpe06}, together with an appropriate semiclassical interpretation of the skein relation.
The present note is an improved version of part of the preprint \cite{Miti-Spera18}.\par
\section{Preliminaries}\par
In this section we concisely review some basic notions related to geometric quantization, tailored to our needs (see e.g. \cite{Wo} for background). First recall that a submanifold $\Lambda$ of a symplectic manifold $(M, \omega)$ is {\it Lagrangian} when the symplectic form 
$\omega$ vanishes thereon and it is of maximal dimension with respect to this property.
If $Q$ is a smooth manifold, then its cotangent space $T^*Q$ is a symplectic manifold
(equipped with a canonical symplectic form). A Lagrangian submanifold $\Lambda \subset T^*Q$ in general position can be described in the following
way (Maslov-H\"ormander {\it Morse family theorem}, see e.g. \cite{Mas,Hor,Gui-Ste}): 
there exists (locally) a smooth function $\phi = \phi(q, a), (q, a) \in Q \times {\mathbb R}^k$ (${\mathbb R}^k$ being a space of auxiliary parameters)
and a submanifold
\begin{equation*}
C_{\phi} = \{ (q, a) \in   Q \times {\mathbb R}^k \, \mid \, d_{a} \phi = 0  \}
\end{equation*}
with $d(d_{a})$ of maximal rank thereon (here $d = d_q + d_{a}$) such that the map
$$
C_{\phi} \ni (q, a) \mapsto   (q,d_q \phi ) \in T^*Q  
$$
is an immersion with image $\Lambda$. 
If the Hessian $H_{a}$ (with respect to the auxiliary variables $a$) is non degenerate, one 
can solve $a = a(q)$ and
define the {\it phase function} $F = F(q) := \phi(q, a(q))$, with $(q, dF(q) ) \in \Lambda$. The covector $dF(q) =: p(q)$
is the momentum at $q$. This fails at the
singular points of the obvious projection $\Lambda \rightarrow Q$, but the singular locus $Z$ (the {\it Maslov cycle})
turns out to be orientable and of
codimension $1$ in $\Lambda$ with $\partial Z$ of codimension $\geq 3$. 
Taking a good open cover $\{ V_i \}_{i \in I}$ of $\Lambda$, and letting $\sigma_i$ be the signature of the Hessian
$H_{a}$  on $V_i \setminus Z$, one readily manufactures the so-called {\it Maslov cocycle}
$\{ h_{ij} = \frac{1}{2} (\sigma_i - \sigma_j) \}$ yielding a class ${\mathcal M} \in H^1(\Lambda,{\mathbb Z})$, dual to the Maslov cycle $Z$, see e.g. \cite[Ch.II, §7]{Gui-Ste}. 
This situation holds for a general symplectic manifold, as a consequence of a result of Weinstein (\cite{Wei}).\par
Now, given a {\it prequantizable} symplectic manifold $(M,\omega)$, i.e. $[\omega] \in H^2(M, {\mathbb Z})$, then, by Weil-Kostant (see e.g. \cite{Wo}), there exists a complex line bundle ${\mathcal L} \to M$ 
(prequantum bundle), equipped with a Hermitian metric and compatible connection $\nabla$ with curvature $\Omega_{\nabla} = -2\pi i \omega$. Since the symplectic 2-form $\omega$ vanishes on any
Lagrangian submanifold $\Lambda \subset M$, any (local) symplectic potential $\vartheta$ 
($d\vartheta = \omega$)
 becomes a closed form thereon, giving a (local) connection form 
pertaining to the restriction of the prequantum connection $\nabla$, denoted by the same symbol. The latter is a {\it flat} connection, and a global covariantly constant section ($\nabla s = 0$) of the (restriction of) the prequantum line bundle - called {\it WKB wave function} - exists if and only if it has trivial holonomy.  
A WKB wave function is subject to sudden phase changes upon crossing the Maslov cycle $Z$ (``passage through a caustic"), governed
by the Maslov cocycle, see e.g. \cite{Gui-Ste,Wo,Mas}. \par
\section{The HOMFLYPT polynomial as a WKB wave function}
 The theory developed in \cite{BeSpe06}, see also \cite{Spe06},  was aimed at placing the construction of the (Abelian) Witten invariant in \cite{Witten89} on firm ground  by avoiding the use of path integrals and it was strongly inspired by the constructions recalled in the preceding section, albeit with modifications dictated by the infinite dimensional environment. We resume it by closely following these papers with appropriate {\it en route} modifications and referring, for the symplectic, hydrodynamical and knot theoretical background, to \cite{Bry,Gui-Ste,Mas,Hor,Arn-Khe,Vas,Kauffman,Ko}.\par
 We shall act within the {\it generalized Brylinski symplectic manifold of oriented mildly singular links} in ${\mathbb R}^3$, $\widehat{Y}$ (allowing a finite number of crossings and finite order tangencies), whose symplectic structure reads, at a generic {\it link} $L$ with components $L_j$, $j=1,\dots n$, represented up to orientation-preserving reparametrizations by smooth embeddings $\gamma_j \in C^{\infty}(S^1, {\mathbb R}^3) \equiv {\mathcal L}{\mathbb R}^3$ with velocities $\dot{\gamma}_j$:
$$
\Omega_L (\cdot , \cdot ) := \sum_{j=1}^n \int_{L_j} \nu (\dot{\gamma}_i, \cdot,\cdot)
$$ 
(with $\nu = dx \wedge dy \wedge dz$ being the standard volume form of ${\mathbb R}^3$).
The manifold consisting of all {\it bona fide} oriented links in ${\mathbb R}^3$ will be denoted by $Y$, and it is clearly non-connected.\par
The volume form $\nu$ can be portrayed as
$$
\nu = dx \wedge dy \wedge dz = d (z\, dx \wedge dy) \equiv d\hat{\theta}
$$
in terms of the (multisymplectic) potential $\hat{\theta}$; the latter transgresses to a (symplectic) potential $\theta$ for $\Omega$, which vanishes identically when restricted on the plane $z=0$. 
The submanifold $\Lambda \subset {\widehat Y}$ consisting of the links on a plane (with indentations keeping track of crossings) is a {\it Lagrangian} one, see \cite{BeSpe06}. \par
Now observe, again following \cite{BeSpe06},  that {\it the links in ${\mathbb R}^3$ can be interpreted as solutions of the {\rm Euler-Lagrange equation} pertaining to a {\rm Chern-Simons Lagrangian} ({\rm helicity}, in hydrodynamical parlance)  with {\rm source} $T_L$ given by the link itself} 
$$
\Phi = \Phi (A, L) := \frac{k}{8\pi} \int_{{\mathbb R}^3} \,A \wedge dA  + \int_L  A \equiv \frac{k}{8\pi} {\mathcal H}(A)  + T_L( A)
$$
with $A$ denoting an Abelian connection (form) on ${\mathbb R}^3$ with curvature $F_A = dA$ - rapidly decaying at infinity to ensure convergence of the integral - and $k$ a non-zero integer or real number.\par
This CS Lagrangian is then taken, as in \cite{BeSpe06}, as a {\it Morse family}, with the {\it auxiliary parameters} (also cf. \cite{Gui-Ste}) given by the Abelian connections.
The Euler-Lagrange equation reads:
$$
{k\over{4\pi}}F_A + T_L = {k\over{4\pi}}dA + T_L = 0
$$
 i.e. one looks  for a connection (viewed as a {\it de Rham current} (\cite{dR})) whose curvature
is concentrated (i.e. $\delta$-like) on $L$. 
The solution can be given in standard vector calculus terms with a so-called Coulomb gauge fixing, ${\rm div} \,{\mathbf A} = 0$, or,
Hodge theoretically, $\delta A = 0$. The (singular) connection $A_L$
 such that $dA_L = T_L$ and $\delta A_L = 0$ can be compactly written in the form
$$
  A_L = -{{4\pi}\over k} {\Delta}^{-1} \delta \, T_L 
$$
where $\Delta$ is the Hodge Laplacian on 1-forms, acting component-wise as the ordinary Laplacian (up to a negative constant)
since we operate in flat space. Existence, in the sense of currents, follows from the 
H\"ormander-\L ojasiewicz theorem, see e.g. \cite{Vla}.
Notice that if we want to insert $A_L$ into $\Phi$ to get a local phase $\phi$ in accordance with the general portrait depicted above (see e.g. \cite{Gui-Ste,BeSpe06}), we are forced to consider {\it ordinary links}. 
Proceeding as in \cite{BeSpe06} (cf. Theorem 3.1 therein) we get, for the local phase $\phi$, the expression (involving the {\it helicity} 
${\mathcal H}(L) = {\mathcal H}(A_L) $ of a {\it framed} link)
$$
\phi (L) = \Phi( A_L, L)  = -\frac{2\pi}{k}  {\mathcal H}(L) \equiv 2\pi \lambda \,\sum_{i,j=1}^n \ell(i,j),
$$
where $\lambda := -1/k$ and
with $\ell(i,j) = \ell(j,i)$ being the {\it Gauss linking number} of components $L_i$ and $L_j$ if $i\neq j$ and where $\ell(j,j)$ is the {\it framing} of $L_j$, equal to $\ell(L_j, L_j^{\prime})$ with $L_j^{\prime}$ being a section of the normal bundle of $L_j$, see e.g. \cite{Rolfsen,Moffatt-Ricca92,BeSpe06,Spe06}. 
A regular projection of a link $L$ onto a plane produces a natural framing called the {\it blackboard framing}, and ${\mathcal H}(L) = w(L)$, the {\it writhe} of $L$.
The helicity can be interpreted, as in \cite{BeSpe06}, as a {\it regularised signature} (cf. Section 2 above) and, as such, it enters the Maslov theory developed therein, cf. Theorem 4.1. 
Since the symplectic potential of Brylinski's form can be taken equal to zero, 
the phase, i.e. the helicity, is (locally) constant, being a topological invariant.  The Lagrangian submanifold 
$\Lambda$ is thence locally given by the graph
$$
(L, d \, {\mathcal H}(L)) = (L, 0)
$$
($d {\mathcal H}(L) = 0$ is the so-called {\it eikonal equation}, see \cite{BeSpe06,Spe06}). We point out that one could equivalently employ, {\it mutatis mutandis}, the Lagrangian submanifold manufactured via the {\it cone construction} of \cite{BeSpe06}.\par
\smallskip
In our context the assumptions of the Weil-Kostant theorem are fulfilled ($\Omega$ is exact) and a covariantly constant section (also called WKB wave function) is just a {\it locally constant function} on $Y$ since, as in \cite{BeSpe06}, we neglect the so-called ``half-form" correction (see e.g. \cite{Wo}).
One must then accommodate passage through the 
corresponding {\it Maslov cycle} ${Z}$, given in our case by the (mildly) {\it singular links possessing exactly one singular point}  causing a sudden jump of writhe (helicity) by $\pm2$ (see again \cite{BeSpe06}, Theorems 4.1 and 5.1), and,
crucially in the link context, we must take into due account the fact that {\it removal of a crossing changes the number of components of a given link and thus places the new link in a {\rm different} connected component of the space $Y$}.  
To this aim, let us consider the following {\it provisional}  wave function (for genuine links $L$)
\begin{equation}
 \psi = \psi (L) := e^{2\pi i \lambda{\mathcal H}(L)} 
\end{equation}
which is a {\it regular isotopy link invariant} (i.e. up to the first Reidemeister move), cf. \cite{BeSpe06}, Theorems 3.1 and 5.1. 
The generic value taken by $\lambda$ (in particular, it can be taken equal to a root of $\pm 1$) avoids trivialities.\\
\begin{wrapfigure}{r}{0.45\textwidth}
	\includegraphics[width=0.43\textwidth]{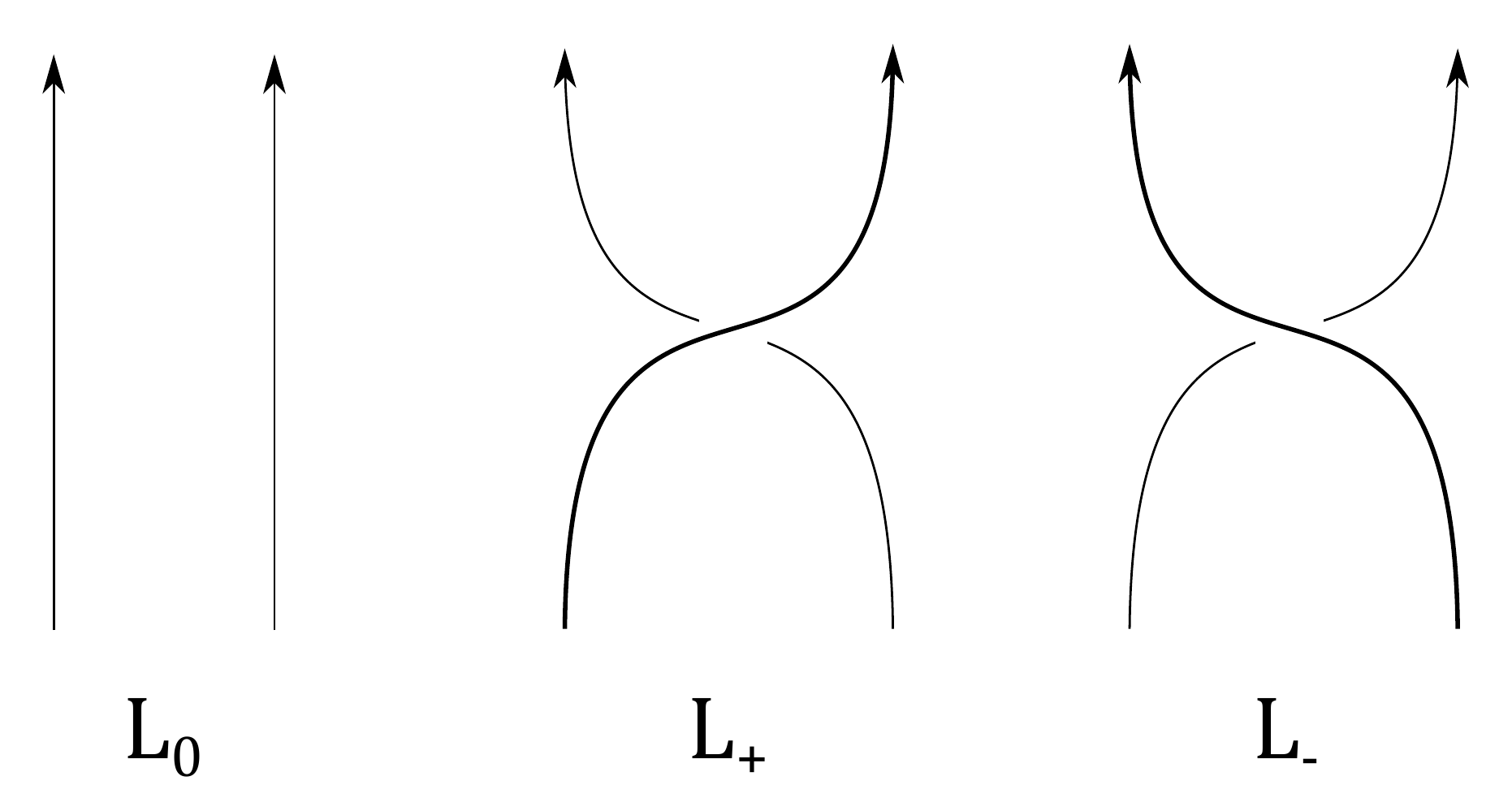}
	\vspace{-3mm}
	\caption{Crossings}
	\label{fig: crossings}
\end{wrapfigure}
Denote, as usual, by $L_+$, $L_-$ and $L_0$ three links (regularly projected onto a plane, $z=0$, say) differing at a single crossing
($(\pm 1)$-crossing, no crossing, respectively), see Figure \ref{fig: crossings}.
Then, inspired by the Liu-Ricca (LR) approach (\cite{Liu-Ricca12,Liu-Ricca15}), let us introduce the ``figures of eight" $E_{\pm}$, that is trivial knots with  $(\pm 1)$-writhe:
${\mathcal H}(E_{\pm}) = 1$.
Starting, for instance, from $L_0$, one can ``add" $E_+$ to the two coherently oriented parallel strands of $L_0$ in such a way that $E_+$ comes with the opposite orientation: a partial cancellation occurs and the net result is 
$L_+$.
\begin{wrapfigure}{l}{0.45\textwidth}
	\includegraphics[width=0.43\textwidth]{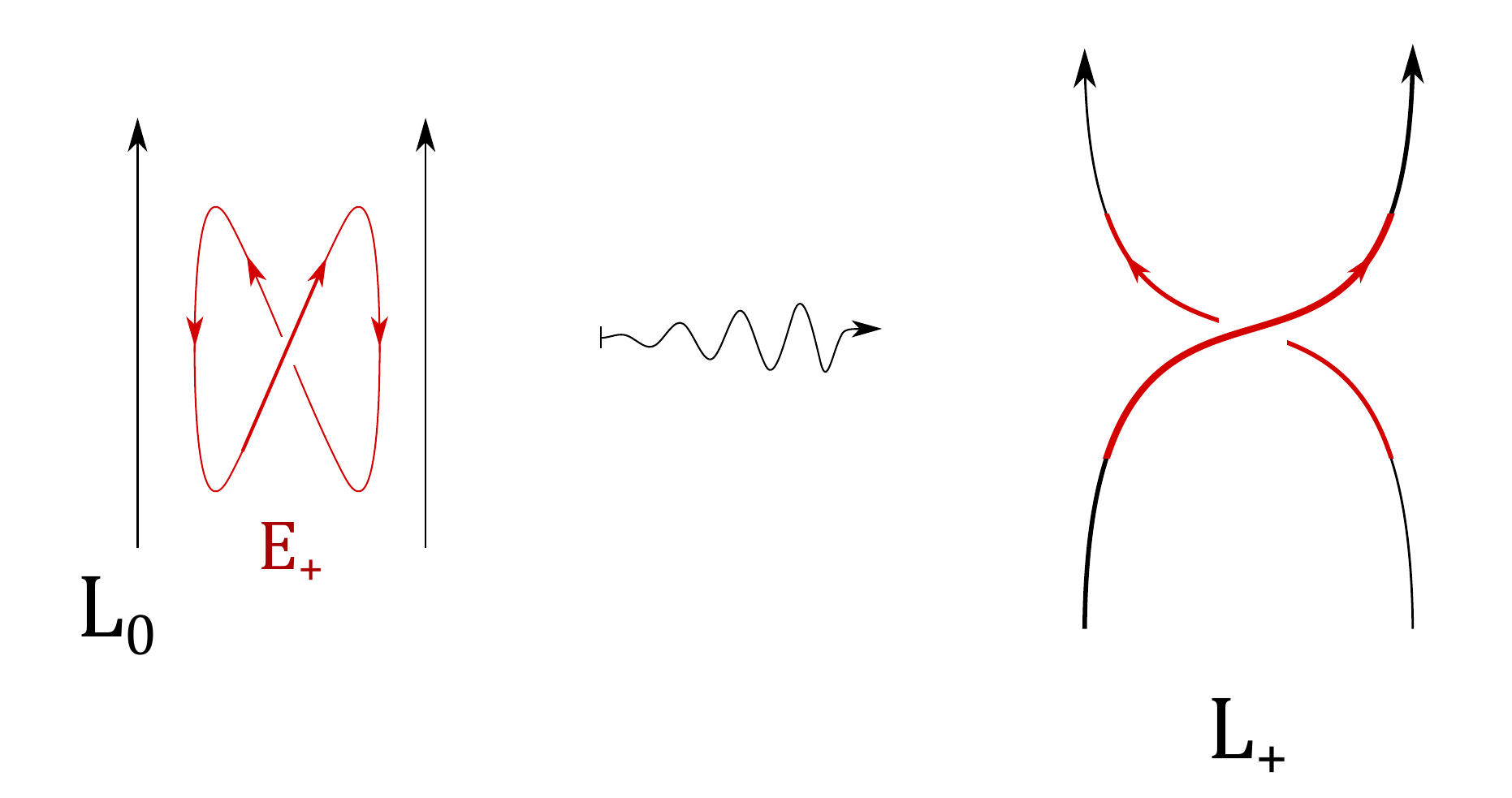}
	\vspace{-3mm}	
	\caption{Surgery via $E_+$}
	\label{fig: surgery}  
	\phantom{.}
\end{wrapfigure}
Conversely, proceeding backwards we can, by adding appropriately an $E_-$, produce $L_0$ from $L_+$
and so on. 
Therefore, addition of $E_{\pm}$ allows one to pass from one local configuration to the other, see Figure \ref{fig: surgery}.\par
Now set: $$
\alpha := e^{2\pi i \lambda {\mathcal H}(E_+)} = e^{2\pi i \lambda}, \quad \qquad  \alpha^{-1} = e^{-2\pi i \lambda} = e^{2\pi i \lambda {\mathcal H}(E_-)}  
$$ 
 so that, trivially, 
 $\psi(L_{\pm}) = \alpha^{\pm 1} \psi(L_{0})$, $ \psi(L_{\pm}) = \alpha^{\pm 2} \psi(L_{\mp})$
 and
 \begin{equation}
 \alpha^{-1} \psi (L_+) - \alpha^{} \psi (L_-) = 0. 
 \end{equation}
Thus we see that $\alpha^{\pm 1}$ arises as the local contribution to the WKB wave function $\psi$ upon addition (surgery) of an eight figure (or ``curl") - which can be applied to a single branch as well (first Reidemeister move) - and $\alpha^{\pm 2}$ as the corresponding contribution upon crossing the Maslov cycle $Z$. \par
{\it We now wish to modify $\psi$ so as to produce a genuine {\rm ambient isotopy} link invariant, keeping the above interpretation}.
 For this purpose, let $\Psi$ be a covariantly constant wave function stemming
from application of the GQ-procedure, normalised in such a way that $\Psi (\bigcirc) = 1$ ($\bigcirc$ being the unknot): as such it is not uniquely determined, since $Y$ is not connected, but
 $\Psi$ can be made to depend naturally on two parameters, the above $\alpha$ and $z$, below. We {\it require} that, upon replacement of $\psi$ by $\Psi$, the modified l.h.s. of (3.2) becomes proportional to $\Psi(L_0)$ (for a suitable constant $z$ which is assumed to be universal, i.e. {\it independent} of the specific link at hand. Consequently, the sought wave function $\Psi$ must satisfy
the {\it skein relation}  (and normalization) for the {\it HOMFLYPT polynomial} $P$ (\cite{Freyd-etal,PT} - here $\alpha^{-1}$  is LR's  $a$)
\begin{equation}
\alpha^{-1} \Psi ({L_+}) - \alpha \Psi ({L_-}) = z \Psi ({L_0}) \, , \qquad \qquad \Psi(\bigcirc) = 1,
\end{equation}
this assuring its existence.
 The trivial wave function $\Psi \equiv 1$ requires $\alpha = 1$ and $z=0$.
 The procedure is still partially {\it ad hoc}, this depending on the non-connectedness of the manifold $Y$.
The skein relation (3.3) can be equivalently  written in the form
$$
 \Psi(L_{-}) =  \alpha^{-2} \Psi(L_{+}) - z \alpha^{-1} \Psi(L_{0})
$$
which tells us that $\Psi(L_{-})$ can be obtained by suitably adding $\Psi(L_{+})$, corrected  by a {\it Maslov type} transition (local surgery via   $\alpha^{-2}$ - one has the same number of link components) and $\Psi(L_{0})$,   {\it corrected by a ``component transition"}   $\alpha^{-1}$ (and multiplied by an extra coefficient $z$).
The latter contribution was absent in   \cite{BeSpe06} since that paper dealt with   {\it knots} only. 
Notice that upon setting $z = \alpha^{-1}  -  \alpha$ and letting $\alpha \to 1$, we get the trivial invariant $\Psi \equiv 1$.\par
\smallskip
{\bf Remarks.} 1. In this way we essentially recover the hydrodynamical portrait of Liu and Ricca \cite{Liu-Ricca12,Liu-Ricca15}, essentially stating that `` $P= t^{\mathcal H}$ " albeit with  a different (and more conceptual) interpretation. In particular, the two parameters used in HOMFLYPT are not quite the same. 
The local surgery operation involves helicity, as in LR, but we portray the latter as a local phase function, governing a component  transition or Maslov, upon squaring it, as in \cite{BeSpe06}.\par
\smallskip
\noindent
2. Passage from $L_{\pm}$ to $L_0$ (and conversely) in ${\widehat Y}$ - abutting, as already remarked, at a change in the number of the link components - involves coalescence
of {\it two} crossings into one and corresponding tangent alignment. This is a sort of ``higher order" contribution beyond the Maslov one.\par
\smallskip
The upshot of the previous discussion is the following:
\vspace{-2mm}
\begin{thm}
The HOMFLYPT polynomial $P = P(\alpha, z)$ can be recovered from the geometric quantization procedure applied to the Brylinski manifold $\widehat{Y}$ and to its Lagrangian subspace $\Lambda$, namely, it coincides (after normalization) with a suitable covariantly constant section $\Psi = \Psi(\alpha,z)$ thereby obtained.
The coefficient $\alpha$ of $P$ is a  phase factor related to the helicity of a standard ``eight-figure" and $z$ comes from accounting for the variation of the number of components of a link.  
\end{thm}
\noindent
{\bf Acknowledgements.}
The authors, both members of the GNSAGA group of INDAM, acknowledge support from Unicatt local D1-funds (ex MIUR 60\% funds). 
They are grateful to Marcello Spera for help with graphics. \par

\end{document}